\documentclass{article}%
\usepackage{amsmath}
\usepackage{amsfonts}
\usepackage{amssymb}
\usepackage{graphicx}%
\newtheorem{theorem}{Theorem}

\newtheorem{lemma}[theorem]{Lemma}

\newtheorem{proposition}[theorem]{Proposition}

\newenvironment{proof}[1][Proof]{\noindent\textbf{#1.} }{\ \rule{0.5em}{0.5em}}
\begin{document}

\title{The Cross Curvature Flow of 3-manifolds with Negative Sectional Curvature}
\author{Bennett Chow\thanks{Partially supported by NSF Grant DMS-0203926.}\\UC San Diego
\and Richard S. Hamilton\\Columbia University}
\date{}
\maketitle

\section{The evolution equation}

We introduce an evolution equation which deforms metrics on 3-manifolds with
sectional curvature of one sign. Given a closed 3-manifold with an initial
metric with negative sectional curvature, we conjecture that this flow will
exist for all time and converge to a hyperbolic metric after a normalization.
We shall establish a monotonicity formula in support of this conjecture. Note
that in contrast to negative sectional curvature, every closed $n\,$-manifold
admits a metric with negative Ricci curvature by the work of Gao and Yau
\cite{Gao and Yau} for $n=3$ and Lohkamp \cite{Lohkamp} for all $n\geq3.$ When
$n\geq4,$ Gromov and Thurston \cite{Gromov and Thurston} have shown that there
exist closed manifolds with arbitrarily pinched negative sectional curvature
which do not admit metrics with constant negative sectional curvature. It is
unknown whether such manifolds admit Einstein metrics. In particular, the
stability result of Ye \cite{Ye} assumes more than curvature pinching. When
$n=3,$ it is an old conjecture, which is also a consequence of the
Geometrization Conjecture, that any closed 3-manifold with negative sectional
curvature admits a hyperbolic metric.

Let $\left(  M,g\right)  $ be a $3$-dimensional Riemannian manifold with
negative sectional curvature. The Einstein tensor is $P_{ij}=R_{ij}-\frac
{1}{2}Rg_{ij}.$ We find it convenient to raise the indices: $P^{ij}%
=g^{ik}g^{j\ell}R_{k\ell}-\frac{1}{2}Rg^{ij}.$ The cross curvature tensor is%
\[
h_{ij}=\left(  \frac{\det P^{k\ell}}{\det g^{k\ell}}\right)  V_{ij}%
\]
where $V_{ij}$ is the inverse of $P^{ij}.$ If the eigenvalues of $P_{ij}$ are
$a=-R_{2323},$ $b=-R_{1313},$ $c=-R_{1212},$ then the eigenvalues of $R_{ij}$
are $-\left(  b+c\right)  ,$ $-\left(  a+c\right)  ,$ $-\left(  a+b\right)  $
and the eigenvalues of $h_{ij}$ are $bc,$ $ac,$ $ab.$ Hence if $\left(
M^{3},g\right)  $ has negative sectional curvature, then both $P_{ij}$ and
$h_{ij}$ are positive definite.

\begin{lemma}
We have the following identities

\begin{enumerate}
\item $\nabla_{i}P^{ij}=0$

\item $h^{ij}\nabla_{i}h_{jk}=\frac{1}{2}h^{ij}\nabla_{k}h_{ij}.$
\end{enumerate}
\end{lemma}

Note that the first identity is the contracted second Bianchi identity. The
above identities imply that the cross curvature tensor is dual to the Ricci
tensor in the following sense.

\begin{lemma}
Let $\left(  M^{n},g\right)  $ be Riemannian manifold.

\begin{enumerate}
\item If the Ricci curvature is positive, then the identity map $\iota:\left(
M,g_{ij}\right)  \rightarrow\left(  M,R_{ij}\right)  $ is harmonic.

\item If $n=3$ and the sectional curvature is negative (or positive), then
$\iota:\left(  M,h_{ij}\right)  \rightarrow\left(  M,g_{ij}\right)  $ is harmonic.
\end{enumerate}
\end{lemma}

With this in mind we define the cross curvature flow (XCF)\footnote{We owe
this nice abbreviation to Ben Andrews.} by%
\[
\frac{\partial}{\partial t}g_{ij}=2h_{ij}%
\]
if the sectional curvature is negative and by $\frac{\partial}{\partial
t}g_{ij}=-2h_{ij}$ if the sectional curvature is positive. Note that the XCF
is fully nonlinear whereas the Ricci flow is quasi-linear.

\section{Short time existence}

Let $\mu_{ijk}$ denote the volume form and raise indices by $\mu^{ijk}%
=g^{ip}g^{jq}g^{kr}\mu_{pqr}.$ Our normalization is such that $\mu_{123}%
=\mu^{123}=1.$ We find it convenient to rewrite $h_{ij}$ as%
\[
h_{ij}=\frac{1}{8}R_{i\ell pq}\mu^{pqk}R_{kjrs}\mu^{rs\ell}.
\]

\begin{lemma}
If $\left(  M,g\right)  $ is a closed $3$-manifold with negative (or positive)
sectional curvature, then for any smooth initial metric a solution to the XCF
exists for a short time.
\end{lemma}

\begin{proof}
We consider the case of negative sectional curvature since the case of
positive sectional curvature is similar. Let $\tilde{g}_{ij}$ denote a
variation of the metric $g_{ij}$ and let tildes also denote the variations of
various curvature tensors. We have%
\[
\tilde{R}_{ijk\ell}=\frac{1}{2}\left(  \frac{\partial^{2}\tilde{g}_{jk}%
}{\partial x^{i}\partial x^{\ell}}+\frac{\partial^{2}\tilde{g}_{i\ell}%
}{\partial x^{j}\partial x^{k}}-\frac{\partial^{2}\tilde{g}_{j\ell}}{\partial
x^{i}\partial x^{k}}-\frac{\partial^{2}\tilde{g}_{ik}}{\partial x^{j}\partial
x^{\ell}}\right)  +\cdots
\]
where the dots denote terms with 1 or less derivatives of the metric. Applying
the equality%
\[
\mu^{pqk}R_{kjrs}\mu^{rs\ell}=-2P^{m\ell}\left(  \delta_{j}^{p}\delta_{m}%
^{q}-\delta_{m}^{p}\delta_{j}^{q}\right)
\]
yields%
\begin{align*}
\tilde{h}_{ij}  &  =-\frac{1}{8}\left(  \frac{\partial^{2}\tilde{g}_{\ell p}%
}{\partial x^{i}\partial x^{q}}+\frac{\partial^{2}\tilde{g}_{iq}}{\partial
x^{\ell}\partial x^{p}}-\frac{\partial^{2}\tilde{g}_{\ell q}}{\partial
x^{i}\partial x^{p}}-\frac{\partial^{2}\tilde{g}_{ip}}{\partial x^{\ell
}\partial x^{q}}\right)  P^{m\ell}\left(  \delta_{j}^{p}\delta_{m}^{q}%
-\delta_{m}^{p}\delta_{j}^{q}\right) \\
&  -\frac{1}{8}\left(  \frac{\partial^{2}\tilde{g}_{\ell p}}{\partial
x^{j}\partial x^{q}}+\frac{\partial^{2}\tilde{g}_{jq}}{\partial x^{\ell
}\partial x^{p}}-\frac{\partial^{2}\tilde{g}_{\ell q}}{\partial x^{j}\partial
x^{p}}-\frac{\partial^{2}\tilde{g}_{jp}}{\partial x^{\ell}\partial x^{q}%
}\right)  P^{m\ell}\left(  \delta_{i}^{p}\delta_{m}^{q}-\delta_{m}^{p}%
\delta_{i}^{q}\right) \\
&  +\cdots.
\end{align*}
The symbol of $E\left(  g\right)  =2h$ is obtained from $\tilde{h}_{ij}$ by
replacing $\frac{\partial}{\partial x^{i}}$ by a cotangent vector $\zeta_{i}$
in the highest (second) order terms%
\[
\sigma DE\left(  g\right)  \left(  \zeta\right)  \tilde{g}_{ij}=-P^{m\ell
}\left(  \zeta_{i}\zeta_{m}\tilde{g}_{\ell j}+\zeta_{\ell}\zeta_{j}\tilde
{g}_{im}-\zeta_{i}\zeta_{j}\tilde{g}_{\ell m}-\zeta_{\ell}\zeta_{m}\tilde
{g}_{ij}\right)  .
\]
Since the sectional curvature is negative, $P^{m\ell}$ is positive and the
eigenvalues of the symbol are nonnegative. One checks that the integrability
condition is $L\left(  h_{ij}\right)  =0,$ where%
\[
L\left(  T\right)  _{k}\doteqdot h^{ij}\nabla_{i}T_{jk}-\frac{1}{2}%
h^{ij}\nabla_{k}T_{ij}.
\]
By Theorem 5.1 of \cite{3d Ric pos}, a solution to the XCF exists for short time.
\end{proof}

\section{Evolution of the Einstein tensor}

We find it convenient to express the Einstein tensor as%
\[
P^{mn}=-\frac{1}{4}\mu^{ijm}\mu^{k\ell n}R_{ijk\ell}.
\]

\begin{lemma}
The evolution of the Einstein tensor is given by%
\[
\frac{\partial}{\partial t}P^{ij}=\nabla_{k}\nabla_{\ell}\left(  P^{k\ell
}P^{ij}-P^{ik}P^{j\ell}\right)  -\det P\,g^{ij}-H\,P^{ij}%
\]
where $H\doteqdot g^{ij}h_{ij}$ and $\det P=\det P^{k\ell}/\det g^{k\ell}.$
\end{lemma}

\begin{proof}
We have%
\begin{align*}
\frac{\partial}{\partial t}R_{ijk\ell}  &  =\nabla_{i}\nabla_{\ell}%
h_{jk}+\nabla_{j}\nabla_{k}h_{i\ell}-\nabla_{i}\nabla_{k}h_{j\ell}-\nabla
_{j}\nabla_{\ell}h_{ik}\\
&  +g^{pq}\left(  R_{ijkp}h_{q\ell}+R_{ijp\ell}h_{qk}\right)  .
\end{align*}
Since the evolution of the volume form is given by $\frac{\partial}{\partial
t}\mu_{ijk}=H\mu_{ijk}$ and $\frac{\partial}{\partial t}\mu^{ijk}=-H\mu
^{ijk},$ we have%
\begin{align*}
\frac{\partial}{\partial t}P^{mn}  &  =-\frac{1}{4}\mu^{ijm}\mu^{k\ell
n}\left(  \nabla_{i}\nabla_{\ell}h_{jk}+\nabla_{j}\nabla_{k}h_{i\ell}%
-\nabla_{i}\nabla_{k}h_{j\ell}-\nabla_{j}\nabla_{\ell}h_{ik}\right) \\
&  -\frac{1}{4}\mu^{ijm}\mu^{k\ell n}g^{pq}\left(  R_{ijkp}h_{q\ell
}+R_{ijp\ell}h_{qk}\right)  -2HP^{mn}\\
&  =\mu^{ijm}\mu^{k\ell n}\nabla_{i}\nabla_{k}h_{j\ell}-\frac{1}{2}\mu
^{ijm}\mu^{k\ell n}g^{pq}R_{ijp\ell}h_{qk}-2HP^{mn}.
\end{align*}
The lemma follows from the identity%
\[
\frac{1}{2}\mu^{ijm}\mu^{k\ell n}g^{pq}R_{ijp\ell}h_{qk}+HP^{mn}=\det
P\,g^{mn}.
\]

\end{proof}

\section{Monotonicity of the volume of the Einstein tensor}

\begin{proposition}
If $\left(  M,g\right)  $ is a $3$-manifold with negative sectional curvature,
then vol$\left(  P_{ij}\right)  $ is nondecreasing under the XCF.
\end{proposition}

This follows from the more general computation

\begin{lemma}
For any $\eta\in\mathbb{R}$%
\begin{align*}
\frac{d}{dt}\int_{M}\left(  \det P\right)  ^{\eta}d\mu &  =\eta\int_{M}\left(
\frac{1}{2}\left\vert T^{ijk}-T^{jik}\right\vert ^{2}-\eta\left\vert
T^{i}\right\vert ^{2}\right)  \left(  \det P\right)  ^{\eta}d\mu\\
&  +\left(  1-2\eta\right)  \int_{M}\left(  \det P\right)  ^{\eta}H\,d\mu
\end{align*}
where $T^{ijk}=P^{i\ell}\nabla_{\ell}P^{jk},$ $T^{i}=V_{jk}T^{ijk}%
=P^{ij}\nabla_{j}\log\det P,$ and the norms are with respect to the metric
$V_{ij}.$
\end{lemma}

Decomposing $T^{ijk}$ into its irreducible components%
\[
T^{ijk}=E^{ijk}-\frac{1}{10}\left(  P^{ij}T^{k}+P^{ik}T^{j}\right)  +\frac
{2}{5}P^{jk}T^{i},
\]
where the coefficients $-\frac{1}{10}$ and $\frac{2}{5}$ are chosen so that
$V_{ij}E^{ijk}=V_{ik}E^{ijk}=V_{jk}E^{ijk},$ we find that%
\[
\left\vert T^{ijk}-T^{jik}\right\vert ^{2}=\left\vert E^{ijk}-E^{jik}%
\right\vert ^{2}+\left\vert T^{i}\right\vert ^{2}.
\]
Taking $\eta=1/2$ in the lemma, we have%
\[
\frac{d}{dt}\int_{M}\left(  \det P\right)  ^{1/2}d\mu=\frac{1}{4}\int
_{M}\left\vert E^{ijk}-E^{jik}\right\vert ^{2}\left(  \det P\right)
^{1/2}d\mu\geq0
\]
and the proposition follows.

\section{Approach to hyperbolic in an integral sense}

We show that an integral measure of the difference of the metric from
hyperbolic is monotone decreasing. Let
\[
J=\int_{M}\left(  \frac{P}{3}-\left(  \det P\right)  ^{1/3}\right)  d\mu
\]
where $P=g_{ij}P^{ij}.$ By the arithmetic-geometric mean inequality, the
integrand is nonnegative, and identically zero if and only if $P_{ij}=\frac
{1}{3}Pg_{ij},$ i.e., $g_{ij}$ has constant curvature.

\begin{theorem}
Under the cross curvature flow $\frac{dJ}{dt}\leq0.$
\end{theorem}

\begin{proof}
We compute $\frac{d}{dt}\int_{M}P\,d\mu=3\int_{M}\det P\,d\mu.$ Combining this
with the previous lemma with $\eta=1/3$ and the decomposition of $T^{ijk}$
into its irreducible components, we have%
\begin{align*}
\frac{dJ}{dt}  &  =-\frac{1}{6}\int_{M}\left(  \left\vert E^{ijk}%
-E^{jik}\right\vert ^{2}+\frac{1}{3}\left\vert T^{i}\right\vert ^{2}\right)
\left(  \det P\right)  ^{1/3}d\mu\\
&  -\int_{M}\left(  \frac{H}{3}-\left(  \det h\right)  ^{1/3}\right)  \left(
\det P\right)  ^{1/3}d\mu
\end{align*}
which is nonpositive ($0$ if and only if $g_{ij}$ has constant negative
sectional curvature).
\end{proof}

\section{A maximum principle estimate}

Let $\square=P^{ij}\nabla_{i}\nabla_{j}=\nabla_{i}\nabla_{j}P^{ij},$ which is
an elliptic operator.

\begin{proposition}%
\[
\frac{\partial}{\partial t}\log\det P=\square\log\det P+\frac{1}{2}\left\vert
T^{ijk}-T^{jik}\right\vert ^{2}-2H.
\]

\end{proposition}

By the maximum principle%
\[
\frac{d}{dt}\min_{M}\log\det P\left(  t\right)  \geq-2\max\left\{  H\left(
x,t\right)  :\det P\left(  x,t\right)  =\min_{M}\det P\left(  t\right)
\right\}  \doteqdot-2X\left(  t\right)  .
\]
Suppose for some $T<\infty$ we have $\inf_{M\times\lbrack0,T)}\det P=0.$ It is
not difficult to show that there then exists a sequence of times
$t_{i}\rightarrow T$ such that $\min_{M}\log\det P\left(  t_{i}\right)
\rightarrow\infty$ and $X\left(  t_{i}\right)  \rightarrow\infty.$ Hence there
exists $x_{i}\in M$ such that $\det P\left(  x_{i},t_{i}\right)  \rightarrow0$
and $H\left(  x_{i},t_{i}\right)  \rightarrow\infty.$ Hence there exists a
sequence of points and times where one of the sectional curvatures tends to
zero and another tends to minus infinity.

\section{Conclusion}

In view of the well-developed theory of Ricci flow \cite{Hamilton} and the
ground breaking work of Perelman \cite{Perelman} on the second author's
program for Ricci flow as an approach to the Thurston Geometrization and
Poincar\'{e} conjectures, it is hopeful that further progress can be made on
the cross curvature flow. Recently Ben Andrews \cite{Ben Andrews} has obtained
new estimates for the cross curvature flow. In the case the universal cover of
the initial 3-manifold is isometrically embedded as a hypersurface in
Euclidean or Minkowski 4-space, the Gauss curvature flow (see \cite{FIrey},
\cite{Tso}, \cite{Chow}, \cite{Ham GCF}, \cite{Andrews GCF} for earlier works
on the GCF) induces the XCF for the metric. In this case Andrews has proved
convergence results. In general, he expects long time existence and
convergence of the XCF to reduce to proving local in time regularity (higher
derivative estimates).

\end{document}